\documentclass{elsart}
\usepackage{flafter,amsmath,amsfonts,cite,psfig,Bold_ams,amssymb}
\usepackage[dvips]{graphicx}
\usepackage{psfrag}
\usepackage{FETI_lam0}
\newtoks\files
\files={abrev2,ref_FETI,ref_DD,ref_EF}

\begin{document}
\begin{frontmatter}
%
%%%%%%%%%%%%%%%%%%%%%%%%%%%%%%%%%%%%%
\title{
On the initial estimate of interface forces in FETI methods
}
%%%%%%%%%%%%%%%%%%%%%%%%%%%%%%%%%%%%%
\author[label1]{Pierre Gosselet}
\author[label1]{Christian Rey}
\author[label2]{Daniel J. Rixen\corauthref{cor1}}
\corauth[cor1]{Corresponding author. E-mail: d.j.rixen@wbmt.tudelft.nl}
\address[label1]{
Laboratoire de Modelisation et Mecanique des Structures (LM2S)\\
FRE 2505 du CNRS\\
8 rue du Capitaine Scott, 75015 Paris, France
}
\address[label2]{
T.U Delft, Faculty of Design, Engineering and Production\\
Engineering  Dynamics\\
Mekelweg 2, 2628 CD Delft, The Netherlands
}
%%%%%%%%%%%%%%%%%%%%%%%%%%%%%%%%%

%\tableofcontents

%\maketitle

%%%%%%%%%%%%%%%%%%%%%%%%%%%%%%%%%%%%%%%%%%%%%%%%%%%%%%%%%%%%%%%%%%%%%
\begin{abstract}
The Balanced Domain Decomposition (BDD) method and the
Finite Element Tearing and Interconnecting (FETI) method are two commonly used
non-overlapping domain decomposition
methods. Due to strong theoretical and numerical similarities,
these two methods are generally considered as being equivalently efficient. However,
for some particular cases, such as for structures with strong
heterogeneities, FETI requires a large number of iterations to compute the solution
compared to BDD.
In this paper, the origin of the bad efficiency of FETI in these particular
cases is traced back to poor initial estimates of the interface stresses.
To improve the estimation of interface forces a novel strategy for splitting
interface forces between neighboring substructures is proposed. The additional
computational cost incurred is not significant. This
yields a new initialization for the FETI method and
restores numerical efficiency which makes FETI comparable to
BDD even for problems where FETI was performing poorly. Various simple test problems are
presented to discuss the efficiency of the proposed strategy and to illustrate the
so-obtained numerical
equivalence between the BDD and FETI solvers.
\end{abstract}
%%%%%%%%%%%%%%%%%%%%%%%%%%%%%%%%%%%%%%%%%%%%%%%%%%%%%%%%%%%%%%%%%%%%%


\end{frontmatter}

\section{Introduction} \label{sec:intro}

Domain Decomposition methods provide a natural framework to solve engineering
problems decomposed into subparts. Such problems can arise for instance because
each subdomain is discretized independently, or because the subdomains represent
different physical domains. Decomposed domains can also be created from an
initial single domain problem in order to make  efficient use of parallel computing
hardware, i.e.\ in order to distribute the computing work on several processors.

Among the domain decomposition methods applied in engineering mechanics to solve
elliptic linear problems, two similar methods have emerged in the last decade as
efficient parallel computing methods: the primal and the dual Schur complement
methods. More specifically, two
procedures have been shown to ensure scalability and robustness:
the Balanced Domain Decomposition (BDD) method
\cite{MANDEL:1993:BAL,LETALLEC:1994:DDM} and the Finite Element Tearing and Interconnecting
(FETI) method \cite{FARHAT:1991:FETI,FARHAT:1994:ADV}.
The BDD method is a primal procedure where
preconditioned conjugate gradient iterations are applied to find the interface
displacements that satisfy the interface equilibrium. In the FETI approach,
it is the interface forces that are searched for iteratively so as to satisfy the
interface displacement compatibility. Therefore FETI is sometimes referred to as a
dual method.

The primal and dual Schur complement methods are based on very
similar concepts (see e.g.\ \cite{RIXEN:2001:PARA} for a
mechanical description). Mathematically, it has been shown that
the preconditioned interface operators for both the BDD and FETI
methods have a condition number bounded by
 \cite{MANDEL:1996:OPT,KLAWONN:2001:BAB}
\begin{equation}
\kappa = O\left(1+\mbox{log}\frac{H}{h}\right)^2
\label{eq:condnum}
\end{equation}
where $H$ and $h$ represent the subdomain and the mesh size respectively.
Hence it is often accepted in the Domain Decomposition community
that using one method instead of the other is a matter of taste
and implementation preferences. However,  when looking closely at
the details and variants of both methods, it becomes rapidly clear
that showing the exact equivalence between the primal BDD and the
dual FETI is not easy, if at all possible
\cite{FRAGAKIS:2002:FETIP1}.

%\cite{FRAGAKIS:2002:FETIP1}, the authors indicate that similar
%conditioning numbers for the BDD and FETI are observed in many cases, but they fall
%short of proving the equivalence.

\begin{figure}[ht]\centering
\begin{minipage}[t]{0.3\textwidth}
\centering
\includegraphics[width=0.9\textwidth,clip=true]{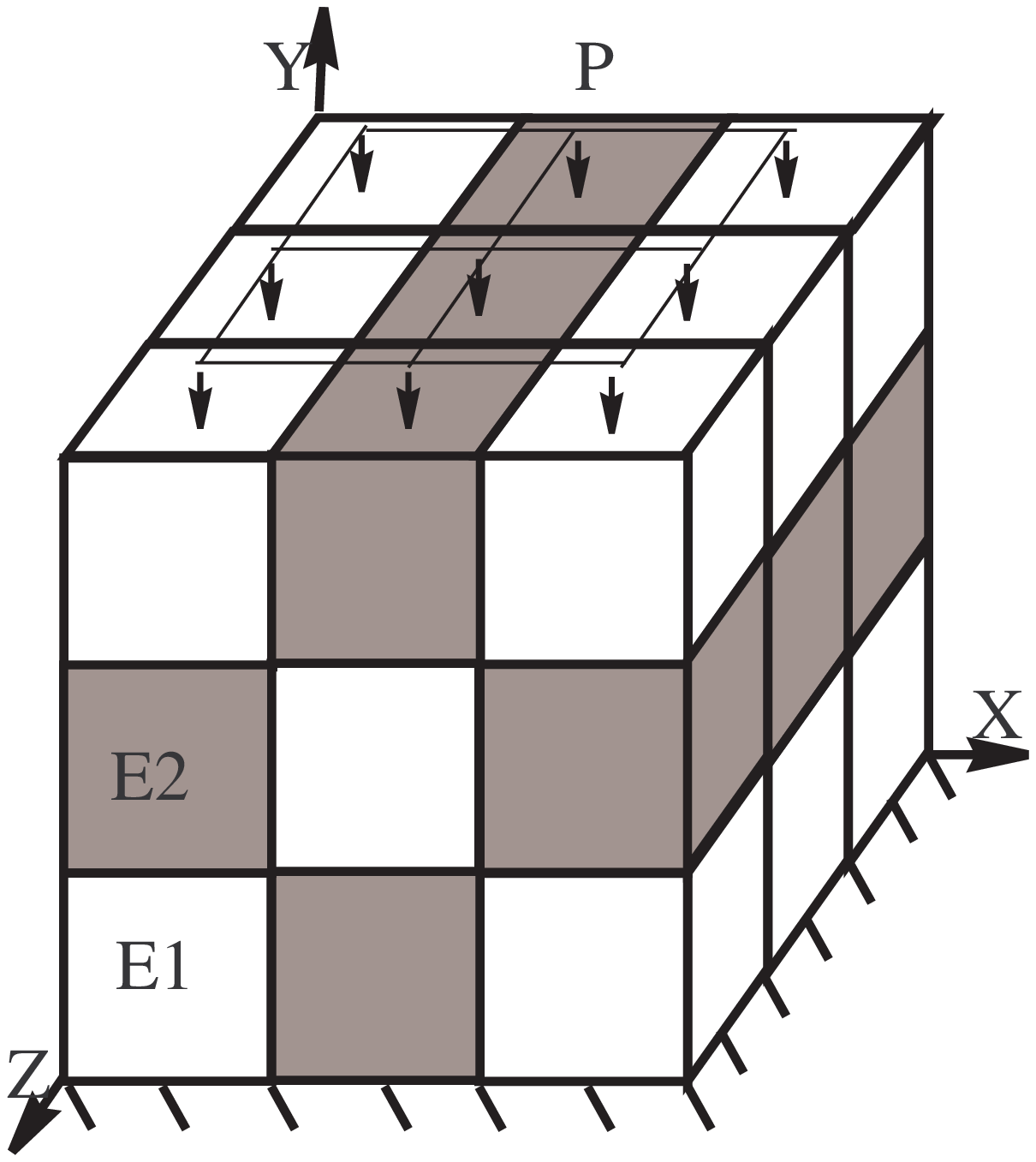}
\caption{Decomposed heterogeneous cube}\label{fig:cube}
\end{minipage}
\begin{minipage}[t]{0.6\textwidth}\vspace*{-3.65cm}
\centering
\psfrag{A}[][][0.75][0]{$\log(\|Primal\ Residual\|)$}
\psfrag{B}[][][0.75][0]{Iteration Number}
\psfrag{C}{}
\psfrag{D}[Tr][][0.5][0]{FETI}
\psfrag{P}[Tr][][0.5][0]{BDD}
\includegraphics[height=0.95\textwidth,angle=-90]{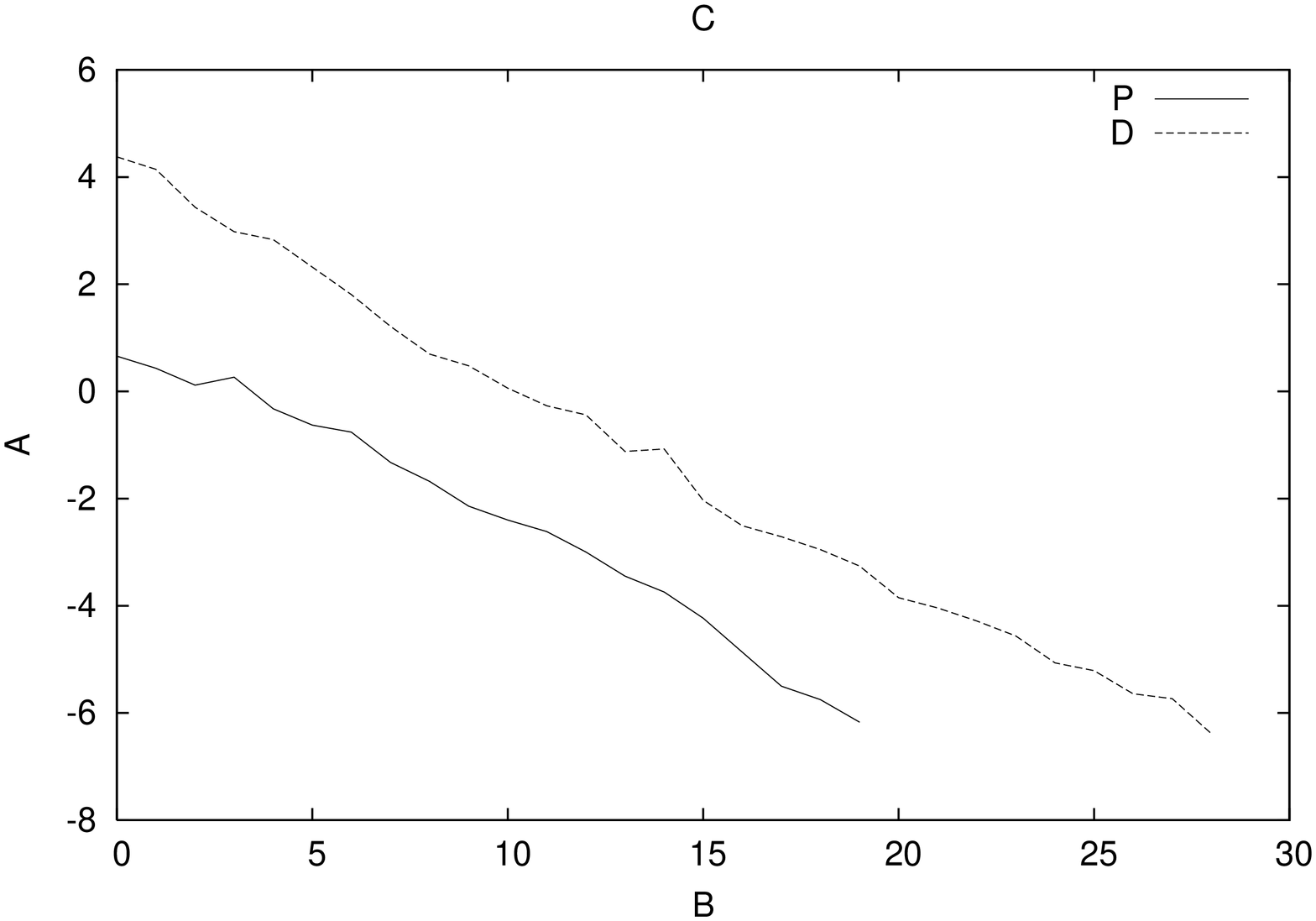}
\caption{Convergence of BDD \& FETI}\label{fig:cubeconv}
\end{minipage}
\end{figure}

Let us then consider the simple example depicted in Figure \ref{fig:cube} of a
highly heterogeneous
elastic cube ($\frac{E1}{E2}=10^5$)
subdivided into $3\times 3\times 3$ subdomains. A uniform pressure is applied on the face
opposite to the clamped side. The structure is discretized
using $Q_2$ hexahedral finite elements with $27$ nodes. The model contains
 $21 000$ degrees of freedom, of which $6 000$  belong to subdomain interfaces.
  In Figure \ref{fig:cubeconv} the convergence
curves of the global equilibrium residual corresponding to iterations of
the BDD and the FETI methods are plotted.
Both methods are equipped with what literature refers to as the best preconditioners
and coarse grids (see section \ref{sec:FETI}). Although the asymptotic convergence of
both procedures is similar, it is observed that the convergence of the FETI method is
less monotonic and that its initial residual is significantly higher.

These results seem to indicate that some specific details such as
the choice of the initial estimates result in possibly significant
differences between the primal and dual algorithms.

In this paper we
revisit the way the interface forces are estimated initially in the
FETI method in an attempt to obtain a convergence comparable
to that of the BDD.

In the next section, we shortly recall the concepts underlying the FETI
solver.
%We also discuss
%the issues of scaling and of coarse grids in a manner different from previous
%publications, in order to give novel insight.
In section \ref{sec:split} we explain that the forces applied on the interface can be
split in different ways. Although the final result is independent of that splitting,
it affects the actual FETI iteration history. We then present an efficient way to define such a
splitting and show how it is related to the construction of the initial iterate of FETI
in section \ref{sec:improved}.
Numerical examples are reported in section \ref{sec:application}
to illustrate the effectiveness of the new initialization strategy.
Finally,
%we will discuss some open issues when
%comparing FETI to the BBD and
we present some conclusions.

%%%%%%%%%%%%%%%%%%%%%%%%%%%%%%%%%%%%%%%%%%%%%%%%%%%%%%%%%%%%%%%%%%%%%
\section{FETI basics}
\label{sec:FETI}

\subsection{The decomposed problem}

Let us consider a domain $\Omega$ subdivided into $N_s$ non-overlapping
subdomains $\Omega^{(s)}$
and assume that we are solving a linear (or linearized) static equilibrium problem on the domain.
The discretized subdomain equilibrium is expressed by
\begin{equation}
\Ks\us = \fs + \gs
\qquad\qquad
s=1,\dots N_s
\label{eq:equis}
\end{equation}
where $\Ks$, $\us$ and $\fs$ are the subdomain stiffness matrices,
displacements and applied forces respectively. $\gs$ are the
connecting forces on the interface between subdomains (thus zero
on the internal degrees of freedom). For the sake of simplicity,
we assume in the following that the meshes are matching
(conforming) on the interface.

The interface forces satisfy an interface equilibrium equation expressing that when assembled
on the interface, the resultant is null (action-reaction):
\begin{equation}
\sum_{s=1}^{N_s} \Ls^T \gs = \zerobold
\label{eq:compgs}
\end{equation}
where $\Ls$ is a Boolean assembly matrix.
The interface connecting forces are such that the interface degrees of freedom are compatible,
namely
\begin{equation}
\sum_{s=1}^{N_s} \Bs\us = \zerobold
\label{eq:compus}
\end{equation}
This relation expresses that for any pair $(u^{(s)},\ u^{(r)})$
of degrees of freedom matching and the interface, $u^{(s)}-u^{(r)}=0$.
$\Bs$ are thus signed Boolean matrices expressing the compatibility constraints on
the interface.

The equilibrium problem of domain $\Omega$ is fully described by the
local equilibrium (\ref{eq:equis}) and by the interface constraints
(\ref{eq:compgs}, \ref{eq:compus}). In block diagonal notations, it can be summarized as
\begin{equation}
\left\{\begin{array}{rcl}
\K \u &=& \f + \g
\\
 \L^T\g &=& \zerobold \\
 \B\u  &=& \zerobold
\end{array}\right.
\label{eq:setufg}
\end{equation}
where $\K$ is the block diagonal matrix of the local operators $\Ks$ and where
\begin{eqnarray}
\u&=&
\left[\begin{array}{c}
\u^{(1)} \\ \vdots \\ \u^{(N_s)}
\end{array}\right]
\qquad
\f=
\left[\begin{array}{c}
\f^{(1)} \\ \vdots \\ \f^{(N_s)}
\end{array}\right]
\qquad
\g=
\left[\begin{array}{c}
\g^{(1)} \\ \vdots \\ \g^{(N_s)}
\end{array}\right]
\nonumber\\
\L^T&=&
\left[\begin{array}{ccc}
\L^{(1)^T} & \cdots & \L^{(N_s)^T}
\end{array}\right]
\nonumber\\
\B&=&
\left[\begin{array}{ccc}
\B^{(1)} & \cdots & \B^{(N_s)}
\end{array}\right]
\end{eqnarray}
Note that in this description, one set of interface displacements and one set
of interface forces are defined per
subdomain.

\subsection{Solvers for decomposed problems}

Solving (\ref{eq:setufg}) can be done in several ways:
\begin{itemize}
\item
Considering (\ref{eq:setufg}) as a constrained equilibrium
problem in terms of $\u$ and $\g$ leads to the three-field formulation
of decomposed domains (see e.g.\ \cite{PARK:1997:ALG, RIXEN:1999:AFETI}).
\item
One can choose to work with a displacement set $\u$ that satisfies
a priori the interface compatibility (\ref{eq:compus}). For that purpose we define a global set $\u_g$ of
degrees of freedom unique on the interface such that
\begin{equation}
\u^{(s)} = \Ls\u_g
\quad\mbox{ or }\quad
\u = \L\u_g
\label{eq:uLug}
\end{equation}
where $\Ls$ is the same assembly Boolean matrix as in
(\ref{eq:compgs}) that extracts subdomain degrees of freedom from
the global set. Stating that $\u^{(s)}$ are obtained from a unique
set is obviously equivalent to stating the interface compatibility
(\ref{eq:compus}) and (\ref{eq:uLug}) thus implies
\begin{equation}
\B\u = \B\L\u_g = \zerobold
\end{equation}
for any global displacement $\u_g$. On the other hand, all compatible displacements
can be written as in (\ref{eq:uLug}).
Hence
\begin{equation}
\L = \mbox{null}(\B)
\label{eq:LnullB}
\end{equation}
In order to illustrate these concepts, we reproduce
the example of \cite{RIXEN:2002:FETI1MPC} in Figure \ref{fig:BL}.
%=================================================================
%  FIGURES
%
%
\begin{figure}
\setlength{\unitlength}{1cm}
\begin{center}
\begin{picture}(4.0,4.5)(0.0,-0.1)
\put(0,10){\psfig{figure=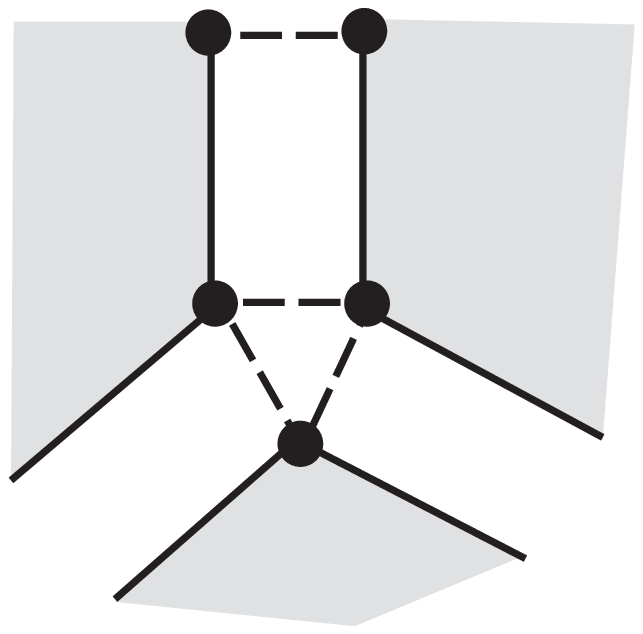,height=3.6cm,width=4cm}}
\put(-0.3,3.6){$\Omega^{(1)}$}
\put(4.1,3.6){$\Omega^{(2)}$}
\put(3.0,-0.0){$\Omega^{(3)}$}
\put(0.5,3.2){$u^{(1)}_1$}
\put(2.6,3.2){$u^{(2)}_1$}
\put(0.6,1.9){$u^{(1)}_2$}
\put(2.7,1.9){$u^{(2)}_2$}
\put(1.7,0.6){$u^{(3)}_1$}
\put(1.7,3.6){$\lambda_1$}
\put(1.7,2.1){$\lambda_2$}
\put(1.2,1.3){$\lambda_3$}
\put(2.3,1.3){$\lambda_4$}
\end{picture}
{\small
\begin{eqnarray}
\underbrace{\hskip3mm
\overbrace{\hskip-3mm
\left[
\begin{array}{cc}
-1 & 0  \\
 0 & -1\\
 0 & -1 \\
 0 & 0
\end{array}\right.
}^{\Bbold^{(1)}}
\overbrace{
\begin{array}{cc}
1 & 0 \\
0 & 1 \\
0 & 0\\
0 & -1
\end{array}
}^{\Bbold^{(2)}}
\overbrace{
\left.\begin{array}{c}
 0\\
 0\\
 1\\
 1
\end{array}
\right]
\hskip-3mm
}^{\Bbold^{(3)}}
\hskip3mm
}_\Bbold
\left[\begin{array}{c}
u^{(1)}_1 \\ u^{(1)}_2 \\
u^{(2)}_1 \\ u^{(2)}_2 \\ u^{(3)}_1
\end{array}\right]
&=&
\zerobold
\nonumber\\
\left[\begin{array}{c}
u^{(1)}_1 \\ u^{(1)}_2 \\
u^{(2)}_1 \\ u^{(2)}_2 \\ u^{(3)}_1
\end{array}\right]
&=&\underbrace{
\left[
\begin{array}{cc}
1 & 0  \\
0 & 1  \\
1 & 0  \\
0 & 1  \\
0 & 1  \\
\end{array}\right]
}_\Lbold
\left[\begin{array}{c}
u_{g_1} \\ u_{g_2}
\end{array}\right]
\nonumber
\end{eqnarray}
}
\caption{Lagrange multipliers and interface compatibility}% (from \cite{RIXEN:2002:FETI1MPC})
\label{fig:BL}
\end{center}
\end{figure}
Introducing (\ref{eq:uLug}) in (\ref{eq:setufg}) yields
\begin{equation}
\left\{\begin{array}{rcl}
\K \L \u_g &=& \f + \g
\\
\L^T\g &=& \zerobold
\end{array}\right.
\label{eq:setufgP}
\end{equation}
This set of equations is at the basis of the primal iterative solution techniques such as the
Primal Schur Complement or the BDD methods \cite{LETALLEC:1991:DDM, MANDEL:1993:BAL}:
iteration schemes are applied to find the displacements $\u_g$ until
the interface equilibrium $\L^T\g=\L^T(\K\L\u_g-\f)=\zerobold$ is satisfied.

\item
One can choose in (\ref{eq:setufg}) a set of interface forces
satisfying a priori the interface equilibrium $\L^T\g = \zerobold$ while keeping redundant
interface degrees of freedom in $\u$.
According to (\ref{eq:LnullB}), such interface forces have the
generic expression
\begin{equation}
\g^{(s)} = - \Bs^T\lambdabold
\quad\mbox{ or }\quad
\g = - \B^T\lambdabold
\end{equation}
$\lambdabold$ are interface forces that act in opposite directions between any pair
of matching degrees of freedom on the interface and are therefore in equilibrium
(see Figure \ref{fig:BL}).
Problem (\ref{eq:setufg}) becomes
\begin{equation}
\left\{\begin{array}{rcl}
\K \u + \B^T\lambdabold &=& \f
\\
\B\u  &=& \zerobold
\end{array}\right.
\label{eq:setufgD}
\end{equation}
Clearly, $\lambdabold$ are the Lagrange multipliers associated to the interface compatibility
constraints.
This form of the decomposed problem is the basis for the dual procedures such as FETI:
iterative algorithms are applied to compute the interface forces $\lambdabold$ such that
the displacements resulting from the subdomain equilibrium are compatible on the interface.
\item
If one chooses interface displacements that are unique on part of the interface while,
on the remainder of the interface, equilibrated connecting forces are defined, one obtains
hybrid primal/dual approaches such as the FETI-DP procedure \cite{FARHAT:2001:FETI_DP}.
\item
If on the entire interface we use displacements and forces that satisfy a linear
combination of the interface equilibrium and compatibility (e.g.\ Robin type of
boundary conditions), one obtains formulations typically used in wave propagation
analysis such as described for instance by Helmholtz equations \cite{DELABOURDONNAYE:1998:FETIH}.
\item
Finally if both the interface equilibrium and compatibility are enforced a priori,
one obtains the fully assembled form
\begin{equation}
\L^T\K\L \u_g = \L^T\f
\label{eq:ass}
\end{equation}

\end{itemize}

%Let us note that in the case where interfaces are non-matching, the relations above remain
%valid, but the interface assembly and constraint matrices ($\L$ and $\B$) are no longer
%Boolean (see for instance \cite{BERNARDI:1990:MORT}).

\subsection{FETI: the dual iterative solver}

In the FETI method \cite{FARHAT:1994:ADV}, the decomposed problem
(\ref{eq:setufgD}) is expressed in terms of interface forces $\lambdabold$: using
the subdomain equilibrium equations to eliminate $\us$,
\begin{equation}
\us = \K^{(s)^+} \left(\fs-\Bs^T\lambdabold\right) - \Rbold^{(s)}\alphabold^{(s)}
\label{eq:us}
\end{equation}
where $\K^{(s)^+}$ is the inverse of $\Ks$ or a generalized inverse if subdomain
$\Omega^{(s)}$ is floating when disconnected from its neighbors. In the latter case,
$\Rbold^{(s)}$ are the associated rigid body modes, their amplitudes $\alphabold^{(s)}$
being determined such that the interface forces are in equilibrium with the applied forces
$\fs$, i.e.\ such that the subdomain equilibrium is well-posed:
\begin{equation}
\Rbold^{(s)^T}\left(\fs - \Bs^T\lambdabold\right) = \zerobold
\label{eq:self}
\end{equation}
Substituting (\ref{eq:us}) into the interface compatibility condition, and taking account
of (\ref{eq:self}), one obtains the dual interface problem
\begin{equation}
\left[\begin{array}{cc}
\FI & \GI \\
\GI^T & \zerobold
\end{array}\right]
\left[\begin{array}{c}
\lambdabold \\
\alphabold
\end{array}\right]
=
\left[\begin{array}{c}
\dbold \\
\ebold
\end{array}\right]
\label{eq:dualinterface}
\end{equation}
where
\begin{eqnarray}
\FI &=& \sum_{s=1}^{N_s} \Bbold^{(s)} \Kbold^{(s)^+} \Bbold^{(s)^T}
= \Bbold \Kbold^+ \Bbold^T
\nonumber\\
\dbold &=& \sum_{s=1}^{N_s} \Bbold^{(s)} \Kbold^{(s)^+} \fbold^{(s)}
= \Bbold \Kbold^+ \fbold
\nonumber\\
\GI &=&
\left[\begin{array}{ccc}
\Bbold^{(1)}\Rbold^{(1)} & \cdots & \Bbold^{(N_s)}\Rbold^{(N_s)}
\end{array}\right]
= \Bbold\Rbold
\nonumber\\
\alphabold &=&
\left[\begin{array}{c}
\alphabold^{(1)} \\ \vdots \\ \alphabold^{(N_s)}
\end{array}\right]
\quad\mbox{and}\quad
\ebold =
\left[\begin{array}{c}
\Rbold^{(1)^T}\fbold^{(1)} \\ \vdots \\ \Rbold^{(N_s)^T}\fbold^{(N_s)}
\end{array}\right]
= \Rbold^T\fbold
\nonumber
\end{eqnarray}
Let use define the projection operator $\P$ such that $\GI^T\P = \zerobold$:
\begin{equation}
\Pbold(\Qbold) = \Ibold - \Qbold\GI\left(\GI^T\Qbold\GI\right)^{-1}\GI^T
\label{eq:defP}
\end{equation}
The choice of the operator $\Qbold$ is discussed later. Introducing the splitting
\begin{eqnarray}
\lambdabold &=& \Pbold(\Q)\bar\lambdabold + \lambdabold_0 \\
\lambdabold_0 &=& \Q\GI\left(\GI^T\Q\GI\right)^{-1}\ebold
\end{eqnarray}
the dual interface problem (\ref{eq:dualinterface})
is equivalent to
\begin{eqnarray}
\P(\Q)^T\FI\P(\Q)\bar\lambdabold = \P(\Q)^T\left(\dbold-\FI\lambdabold_0\right)
\label{eq:dualPQ}
\end{eqnarray}

The FETI method consists in  preconditioned conjugate gradient iterations on
the dual interface problem (\ref{eq:dualPQ}) \cite{FARHAT:1994:ADV}.
Applying $\FI$ to an iterate is a naturally parallel operation. The projection steps
however require solving a coarse global problem.
The preconditioning operators and
the possible choices for $\Q$ are summarized next.

\subsection{FETI preconditioners and coarse grid space}
\label{subsec:precond}

The forces $\B^T\lambdabold$ exist only on  the interface degrees of freedom. Using
a subscript $_b$ and $_i$ for interface boundary and internal degrees of freedom
respectively, one can re-write the decomposed problem (\ref{eq:setufgD}) by condensing
out $\u_i$:
\begin{equation}
\left\{\begin{array}{rcl} \S \u_b + \Bb^T\lambdabold &=& \f_b^*
\\
\Bb\u_b  &=& \zerobold
\end{array}\right.
\label{eq:setufgDS}
\end{equation}
where $\Bb$ is a submatrix of $\B$ associated to boundary
displacements only and
\begin{eqnarray}
\S &=&
\left[\begin{array}{ccc}
\quad\qquad\ddots \\
&\hspace{-10mm}
\K^{(s)}_{bb}-\K^{(s)}_{bi}\K^{(s)^{-1}}_{ii}\K^{(s)}_{ib}
\hspace{-10mm} \\
&&\ddots\quad\qquad
\end{array}\right]
\\
\f_b^* &=&
\left[\begin{array}{c}
\vdots \\
\f^{(s)}_{b}-\K^{(s)}_{bi}\K^{(s)^{-1}}_{ii}\f^{(s)}_{i} \\
\vdots
\end{array}\right]
\end{eqnarray}
$\S$ is the block diagonal matrix of the local operators
statically condensed on the interface (also known as Schur
complements). The dual interface problem  (\ref{eq:dualinterface})
or (\ref{eq:dualPQ}) can thus be written in equivalent forms by
replacing $\K$, $\f$ and $\B$ by $\S$, $\f_b^*$ and $\Bb$
respectively.

The preconditioner in FETI for which the optimal conditioning number (\ref{eq:condnum})
holds is the Dirichlet preconditioner:
\begin{equation}
\tilde\Fbold_{I_D}^{-1}= \tilde\Bb\S\tilde\Bb^T \label{eq:DirPrec}
\end{equation}
The operator $\tilde\Bb$ is similar to $\Bb$ but includes a
scaling with respect to interface multiplicity or relative
interface stiffness \cite{RIXEN:1998:SUPERL}. The scaling is
implemented as a simple pre- and post-processing in the
preconditioning step. Mathematically, its general expression is
\cite{KLAWONN:2001:BAB, RIXEN:2002:FETI1MPC}
\begin{equation}
\tilde\Bb = \left(\Bb\Abold\Bb^T\right)^{+}\Bb\Abold
\end{equation}
where the $^+$ superscript denotes an inverse or a pseudo-inverse in case
redundant compatibility constraints are present. It can be shown that the multiplicity
scaling procedure corresponds to $\Abold=\Ibold$, whereas the stiffness
scaling (also known as super-lumped scaling)
corresponds to $\Abold=\mbox{diag}(\K_{bb})^{-1}$.

The operator $\Qbold$ in the coarse grid projector (\ref{eq:defP}) defines
the interface forces $\Qbold\GI$ associated to rigid mode displacements of
subdomains. It can be chosen as $\Qbold=\Ibold$, but for many engineering problems
it should be taken as the preconditioner, namely
\begin{equation}
\Qbold=\tilde\Fbold_{I_D}^{-1}
\label{eq:Q1}
\end{equation}
 or
at least as a degenerated form of the preconditioner such as
\begin{equation}
\Qbold=\left(\Bb\ \mbox{diag}(\K_{bb})^{-1}\Bb^T\right)^{+}
\label{eq:Q2}
\end{equation}

Further discussion
on preconditioning, scaling and on  the choice of $\Qbold$ can be found in
\cite{FARHAT:1994:ADV,RIXEN:1998:SUPERL,Bhardwaj:1999:PQ_LRV,RIXEN:2002:FETI1MPC}.

%%%%%%%%%%%%%%%%%%%%%%%%%%%%%%%%%%%%%%%%%%%%%%%%%%%%%%%%%%%%%%%%%%%%%
\section{Splitting forces on an interface}
\label{sec:split}

The definition of the decomposed hybrid problem (\ref{eq:setufgD}) is
not unique: any interface force satisfying the interface equilibrium can be added to the
system without changing the final solution $\u$. Indeed, replacing (\ref{eq:setufgD})
by
\begin{equation}
\left\{\begin{array}{rcl}
\K \u + \B^T\tilde\lambdabold &=& \tilde\f \  =\  \f+\B^T\mubold
\\
\B\u  &=& \zerobold
\end{array}\right.
\label{eq:setufgDmu}
\end{equation}
for any $\mubold$ will yield the same solution $\u$, the Lagrange multiplier being
then such that $\tilde\lambdabold=\lambdabold+\mubold$. This can also be observed from
the fact that (\ref{eq:setufgD}) and (\ref{eq:setufgDmu}) have the same assembled forces
since $\L^T\tilde\f=\L^T\f$. Mechanically speaking, it means that adding forces
on one side of the interface boundary and subtracting it on the other side yields the same
solution in terms of displacement although it modifies the internal forces on the interface.

\subsection{Force splitting commonly used}

From the discussion above, it is clear that one can
choose the splitting of interface forces on the interface
boundaries\footnote{In the same manner, it was
noted in \cite{RIXEN:2002:FETI1MPC} that coefficients of additional constraints can
also be arbitrarily split on the interface, which led to the construction of efficient
preconditioners for problems with multipoint constraints.}.
Commonly, interface forces are split between the
connecting subdomains proportionally to
subdomain relative stiffness, in a way consistent with the scaling in the preconditioner.
Mathematically speaking, this is equivalent to choosing
\begin{equation}
\tilde\f = \diag(\K)\L \left(\L^T\diag(\K)\L\right)^{-1}\left(\L^T\f\right)
\label{eq:splitf}
\end{equation}
where $\left(\L^T\f\right)$ are the applied forces given in an initially assembled
problem or the assembled forces corresponding to a different force splitting.
$(\L^T\diag(\K)\L)$ is the assembled  diagonal stiffness
on the interface.

%In order to show that the choice (\ref{eq:splitf}) corresponds to a particular
%case of the general form (\ref{eq:setufgDmu}), let us transform it as
%%
%\begin{equation}
%\tilde\f = \f +
%\left(\diag(\K)\L \left(\L^T\diag(\K)\L\right)^{-1}\L^T-\Ibold\right) \f
%\label{eq:interm1}
%\end{equation}
%%
Let us note that for any symmetric positive definite matrix $\A$,
the following relation holds: \footnote{ This relation expresses
the complementarity of the primal and dual scaling and was already
suggested in \cite{FRAGAKIS:2002:FETIP1}. Its proof can be
obtained by noting that $\L^T(\ref{eq:complDPscal})$ and
$\B\A^{-1}(\ref{eq:complDPscal})$ are trivially satisfied. Owing
to the fact that $\A$ is symmetric positive definite and because
the image of $\L$ corresponds to the nullspace of $\B$, the final
result is obtained.}
\begin{equation}
\A\L \left(\L^T\A\L\right)^{-1}\L^T
+
\B^T\left(\B\A^{-1}\B^T\right)^{+}\B\A^{-1} = \Ibold
\label{eq:complDPscal}
\end{equation}
From (\ref{eq:complDPscal}),
one deduces that the force splitting  (\ref{eq:splitf}) is equivalent to
\begin{equation}
\tilde\f = \f - \B^T \ \left(\B\diag(\K)^{-1}\B^T\right)^{+}\B\diag(\K)^{-1} \f
\label{eq:splitfB}
\end{equation}
This last relation shows that the common
splitting technique (\ref{eq:splitf})
corresponds to adding a particular set of equilibrated interface forces as
described by (\ref{eq:setufgDmu}).

\subsection{Splitting statically condensed forces}

In section \ref{subsec:precond}, we indicated that the hybrid decomposed
problem can be set in an equivalent form
condensed on the interface (\ref{eq:setufgDS}).
Following the same splitting procedure for the force as in (\ref{eq:splitf}),
one would then construct the problem
\begin{equation}
\left\{\begin{array}{rcl} \S \u_b + \Bb^T\tilde\lambdabold &=&
\tilde\f^*_b
\\
\Bb\u_b  &=& \zerobold
\end{array}\right.
\label{eq:setufgDbmu}
\end{equation}
where
\begin{equation}
\tilde\f^*_b = \diag(\K_{bb})\L_b
\left(\L_b^T\diag(\K_{bb})\L_b\right)^{-1}(\L_b^T\f^*_b)
\label{eq:splitfcond}
\end{equation}
and where $(\L_b^T\f^*_b)$ are the applied forces statically
condensed and assembled on the interface. Following a similar
discussion as in the previous section, it is straightforward to
show that such a splitting is equivalent to
\begin{equation}
\tilde\f^*_b = \f^*_b - \Bb^T \
\left(\Bb\diag(\K_{bb})^{-1}\Bb^T\right)^{+}\Bb\diag(\K_{bb})^{-1}
\f^*_b \label{eq:splitfcond2}
\end{equation}
indicating again that the splitting corresponds to adding a particular
set of equilibrated interface forces.

Let us observe that, in the non-condensed format, (\ref{eq:setufgDbmu}) is identical to
\begin{eqnarray}
%\left[\begin{array}{cc}
%\K_{ii} & \K_{ib} \\
%\K_{bi} & \K_{bb}
%\end{array}\right]
%\left[\begin{array}{c}
% \u_i \\ \u_b
%\end{array}\right]
\K\u
 +
%\left[\begin{array}{c}
% \zerobold \\ \Bb^T
%\end{array}\right]
\B^T
\lambdabold &=&
\left[\begin{array}{c}
 \f_i \\ \tilde\f_b^*+\K_{bi}\K^{-1}_{ii}\f_{i}
\end{array}\right]
\label{eq:setufgDDbmunc}\\
&=&
\f -
\B^T \ \left(\B\diag(\K)^{-1}\B^T\right)^{+}\B\diag(\K)^{-1}
\left[\begin{array}{c}
\zerobold \\
\f^*_b
\end{array}\right]
\nonumber
\end{eqnarray}
When comparing it with (\ref{eq:splitfB}), it is clear that
splitting the condensed forces is not equivalent to splitting $\f$
as commonly done. Although the displacement field  obtained is
obviously the same, the Lagrange multipliers searched for during
the FETI iterations will be different. The question thus arises:
what splitting should be used and how does it affect the FETI
iterations on the dual interface problem? Next section provides a
theoretical analysis of the interest of the splitting while the
section after provides related numerical assessments.

%%%%%%%%%%%%%%%%%%%%%%%%%%%%%%%%%%%%%%%%%%%%%%%%%%%%%%%%%%%%%%%%%%%%%
\section{An improved initial estimate of the Lagrange multipliers}
\label{sec:improved}

Let us consider again the hybrid decomposed problem (\ref{eq:setufgD}). In
order to construct an estimate for the interface Lagrange multipliers,
let us assume
that the internal degrees of freedom $\u_i$
satisfy the local equilibrium while $\u_b$ on
the interface boundary have a zero estimate. This is exactly
the assumption underlying the initialization of the
primal Schur complement iteration schemes. The subdomain equilibrium
then writes
\begin{equation}
\Bb^T\lambdabold  \simeq \f^*_b
\end{equation}
which cannot be exactly satisfied (unless $\u_b=0$ corresponds to
the solution). Hence, let us decompose the statically condensed
force $\f^*_b$ into a component satisfying the interface
equilibrium (i.e.\ belonging to the image of $\Bb^T$) and a
remainder:
\begin{equation}
\Bb^T\lambdabold  \simeq \f^*_b = \Bb^T\gammabold +
\mbox{\textsf{remainder}}
\end{equation}
where
\begin{equation}
\gammabold = \left(\Bb\D\Bb^T\right)^{+}\Bb\D\f^*_b
\label{eq:gamma}
\end{equation}
for any non-singular and symmetric matrix $\D$. The remainder is
then orthogonal to $\D\Bb^T$. If we choose as initial estimate
\begin{equation}
\lambdabold_{00} = \gammabold
\end{equation}
then the initial equilibrium residual
$\Bb^T\lambdabold_{00}-\f^*_b=\mbox{\textsf{remainder}}$ is
minimum in the $\D$-norm. Setting $\D$ to be $\diag(\K_{bb})^{-1}$
\begin{equation}
\lambdabold_{00} =
\left(\Bb\diag(\K_{bb})^{-1}\Bb^T\right)^{+}\Bb\diag(\K_{bb})^{-1}\f^*_b
\label{eq:l00}
\end{equation}
can be easily computed by a scaling of the interface condensed forces. Observing
that $\diag(\K_{bb})^{-1}$ is an approximation of
$\S^+$, $\lambdabold_0$ minimizes a norm having the meaning of an energy.

We thus propose to start the FETI iterations with an initial estimate obtained
by rendering $\lambdabold_{00}$ admissible, i.e. such that $\GI^T\lambdabold_0=\ebold$.
Using the notations introduced in (\ref{eq:defP}),
\begin{equation}
\lambdabold_0 = \Pbold(\Qbold) \lambdabold_{00} +  \Qbold\GI\left(\GI^T\Qbold\GI\right)^{-1}\ebold
\label{eq:l0}
\end{equation}
To be consistent with the choice $\Dbold=\diag(\K_{bb})^{-1}$, one should take $\Qbold$
as in (\ref{eq:Q1}) or (\ref{eq:Q2}).

\subsection*{Remarks}

\begin{itemize}
\item
If $\Dbold=\Sbold^+$ in (\ref{eq:gamma}),
$\lambdabold_{00} = \FI^{+}\dbold$. Furthermore, one would have
$\Qbold=\FI^{+}$ and
\begin{eqnarray}
\lambdabold_0 &=& \Pbold(\Q)\FI^{+}\dbold +  \FI^{+}\GI\left(\GI^T\FI^{+}\GI\right)^{-1}\ebold
\nonumber\\
&=& \FI^{+}\left(\dbold - \GI\alphabold\right) \\
&\mbox{with}&\alphabold = \left(\GI^T\FI^{+}\GI\right)^{-1} (\GI^T\FI^{+}\dbold-\ebold)
\nonumber
\end{eqnarray}
Hence the new choice of initial estimate would yield the exact solution.
\item
If we choose to split the interface forces such as described in (\ref{eq:splitfcond})
(or its equivalent form \ref{eq:splitfcond2}),
(\ref{eq:l00}) yields $\lambdabold_{00} = \zerobold$ so that the initial iterate
(\ref{eq:l0}) would correspond to the standard FETI starting procedure.
Hence \emph{it is equivalent to apply standard FETI iterations
when splitting condensed forces on the interface as in (\ref{eq:splitfcond})
or to use any decomposed
force vector $\f$ together with the starting procedure (\ref{eq:l0}).}
\item
The construction of $\lambdabold_{00}$ in (\ref{eq:l00}) has exactly the same mechanical
interpretation as the scaling of interface forces in the preconditioner (see
\cite{RIXEN:2002:FETI1MPC}).
The cost incurred by (\ref{eq:splitfcond}) or (\ref{eq:l0}), corresponds to a preconditioning
step, thus less than one half of a FETI iteration. Obviously, if the lumped preconditioner
is applied, the standard splitting (\ref{eq:splitf}) should be considered.
\end{itemize}

%%%%%%%%%%%%%%%%%%%%%%%%%%%%%%%%%%%%%%%%%%%%%%%%%%%%%%%%%%%%%%%%%%%%%
\section{Numerical assessment}
\label{sec:application}

In order to assess the performance of the new estimate for the Lagrange multipliers, we consider the solution to various problems by FETI and BDD.
In all cases, convergence is monitored through the evaluation of the global residual:
\begin{equation}
% \text{convergence} \Leftrightarrow
\frac{\| \K_g \u_g - \f_g \| }{ \| \f_g \| }  \leqslant \varepsilon
\end{equation}
where $\K_g$, $\u_g$ and $\f_g$ are global assembled stiffness
matrix, displacement field and forces (see (\ref{eq:ass})).
All results are obtained with the Dirichlet preconditioner equipped with
the super-lumped scaling, i.e.\ $\Abold = \diag(\K_{bb})^{-1}$.
The stopping criterion $\varepsilon$ is set to $10^{-6}$.
The FETI method is tested for its different projectors:
we note $\P(\tilde{\Fbold}^{-1}_{I_D})$ the Dirichlet projector,
$\P(\Wbold)$ the superlumped projector corresponding to (\ref{eq:Q2}) and $\P(\Ibold)$ the identity projector.

\begin{figure}[ht]\centering
\psfrag{A}[][][0.75][0]{$\log(\|Primal\ Residual\|)$}
\psfrag{B}[][][0.75][0]{Iteration Number}
%\psfrag{C}{}
\psfrag{D}[Tr][][0.5][0]{FETI $P(\tilde{\Fbold}^{-1}_{I_D})$}
\psfrag{W}[Tr][][0.5][0]{FETI $\P(\Wbold)$}
\psfrag{I}[Tr][][0.5][0]{FETI $\P(\Ibold)$}
\psfrag{P}[Tr][][0.5][0]{BDD}
\includegraphics[height=0.72\textwidth,angle=-90]{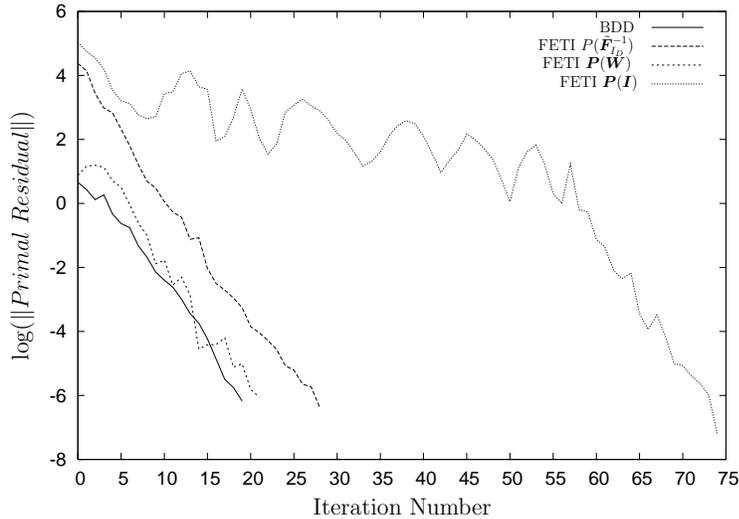}
\caption{Cube with checkerboard heterogeneity (fig.\ \ref{fig:cube}):
 convergence of the primal residual for FETI with standard force splitting}
\label{fig:PRclassical}
\end{figure}

\paragraph{Cube with checkerboard heterogeneity: }
Let us first assess the strategies described above on the problem described in the
introduction (see fig.\ \ref{fig:cube}).
Figure \ref{fig:PRclassical} presents the convergence
of the primal residual through the conjugate gradient iterations for the BDD
and for the FETI approaches with standard force splitting.
The primal approach clearly exhibit better results.
Also, due to the strong heterogeneity in the structure and along the interfaces,
the projector $\P(\Ibold)$ yields very poor convergence and is thus not suitable.
It is observed that,
although the Dirichlet projector yields a better convergence rate, the
superlumped projector $\P(\Wbold)$ has a significantly lower initial residual and hence reaches
convergence faster.
The initialization associated to $\P(\Wbold)$ leads to a starting residual
$10^3$ times smaller,
however its behaviour during the solution process is less regular.

\begin{figure}[ht]\centering
\psfrag{A}[][][0.75][0]{$\log(\|Primal\ Residual\|)$}
\psfrag{B}[][][0.75][0]{Iteration Number}
%\psfrag{C}{}
\psfrag{F}[Tr][][0.5][0]{improved initialization FETI $P(\tilde{\Fbold}^{-1}_{I_D})$}
\psfrag{D}[Tr][][0.5][0]{FETI $P(\tilde{\Fbold}^{-1}_{I_D})$}
\psfrag{P}[Tr][][0.5][0]{BDD}
\includegraphics[height=0.72\textwidth,angle=-90]{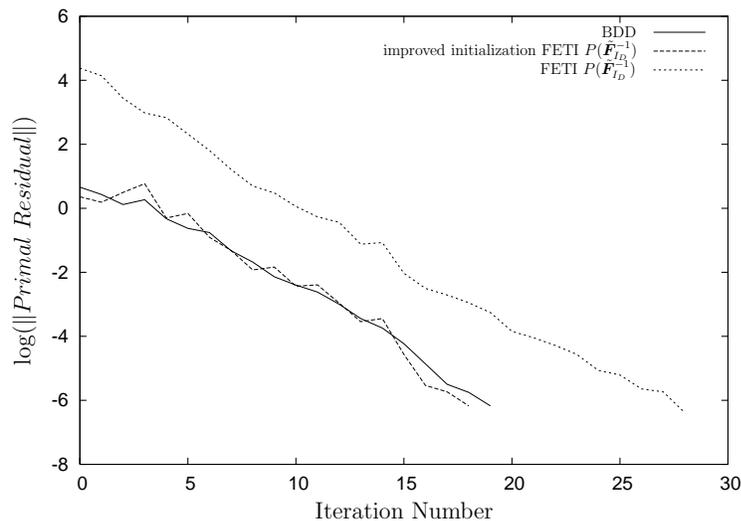}
\caption{Cube with checkerboard heterogeneity (fig. \ref{fig:cube}):
 convergence of the primal residual for FETI with new initialization}
\label{fig:PRdualPF}
\end{figure}

Figure \ref{fig:PRdualPF} presents the convergence history for the BDD method
and for the FETI method equipped with the Dirichlet projector and with the new initialization.
For comparison, the convergence obtained with the standard splitting is shown once more.
Note that the new initialization not only provides a better initial estimate
(lower initial residual) but also leads to similar convergence rates as for
the standard splitting. Therefore the convergence history of FETI becomes
very similar to the convergence of the primal BDD approach.
The total number of iterations is about the same as for the primal BDD, although
the convergence of FETI is less monotonic.

\begin{table}[ht]\centering
\begin{tabular}{|c|c|c|c|}\hline \multicolumn{2}{|c|}{Approach} & Number of  & Initial Residual   \\
\multicolumn{2}{|c|}{ } &  iterations & (Log)  \\ \hline  \hline
\multicolumn{2}{|c|}{BDD} & $19$ & $0.659$ \\ \hline  \hline
FETI       & No Splitting             & $28$  & $4.428$ \\ \cline{2-4}
Dirichlet  & Classical Splitting  & $28$ & $4.377$ \\ \cline{2-4}
 $P(\tilde{F}^{-1}_{I_D})$         & New Initialization    & $18$ & $0.359$ \\ \hline  \hline
FETI        & No Splitting             & $21$ & $0.873$ \\ \cline{2-4}
Dirichlet  & Classical Splitting  & $21$ & $0.872$ \\ \cline{2-4}
$P(W)$        & New Initialization   & $20$ & $0.868$ \\ \hline  \hline
FETI        & No Splitting            & $74$ & $5.029$ \\ \cline{2-4}
Dirichlet  & Classical Splitting  & $74$ & $5.026$ \\ \cline{2-4}
$P(I) $         & New Initialization    & $73$ & $5.016$ \\ \hline
\end{tabular}\caption{Performance results for the cube of fig. \ref{fig:cube}}\label{tab:prob1}
\end{table}

Table \ref{tab:prob1} summarizes the performance results
of the various available strategies.
%As explained in section \ref{sec:split} and \ref{sec:improved}
%the new initialization is equivalent to
%the splitting of condensed forces, we then compare the performance results of the new initialization to the performance results of the no splitting approach and
%of the classical splitting of applied forces appoach.
To investigate the efficiency of the new initialization (or force splitting) procedure,
we compare the number of iterations to
achieve convergence and the norm of the initial primal residual.
The new initialization yields significant improvements
for the Dirichlet projector ($35\%$ less iterations), but
only slightly affects the convergence when other projectors are used.
%In this example, the classical splitting is almost useless. \\

\begin{figure}[ht]\centering
\begin{minipage}[t]{0.4\textwidth}
\centering
\includegraphics[height=0.9\textwidth,clip=true]{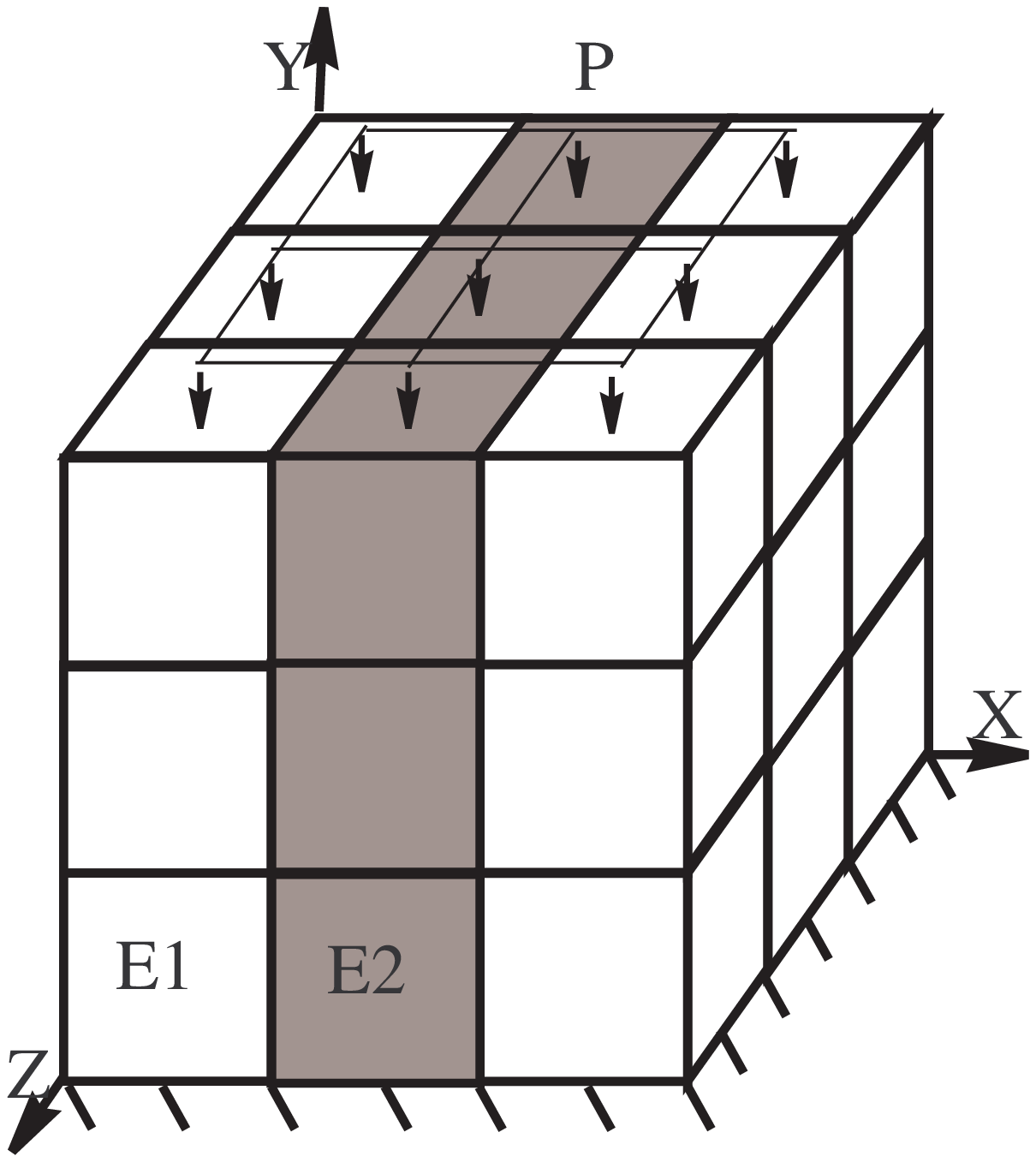}
\caption{Cube with heterogeneous layers}\label{fig:cubeh3}
\end{minipage}
\begin{minipage}[t]{0.4\textwidth}
\centering
\includegraphics[height=0.9\textwidth,clip=true]{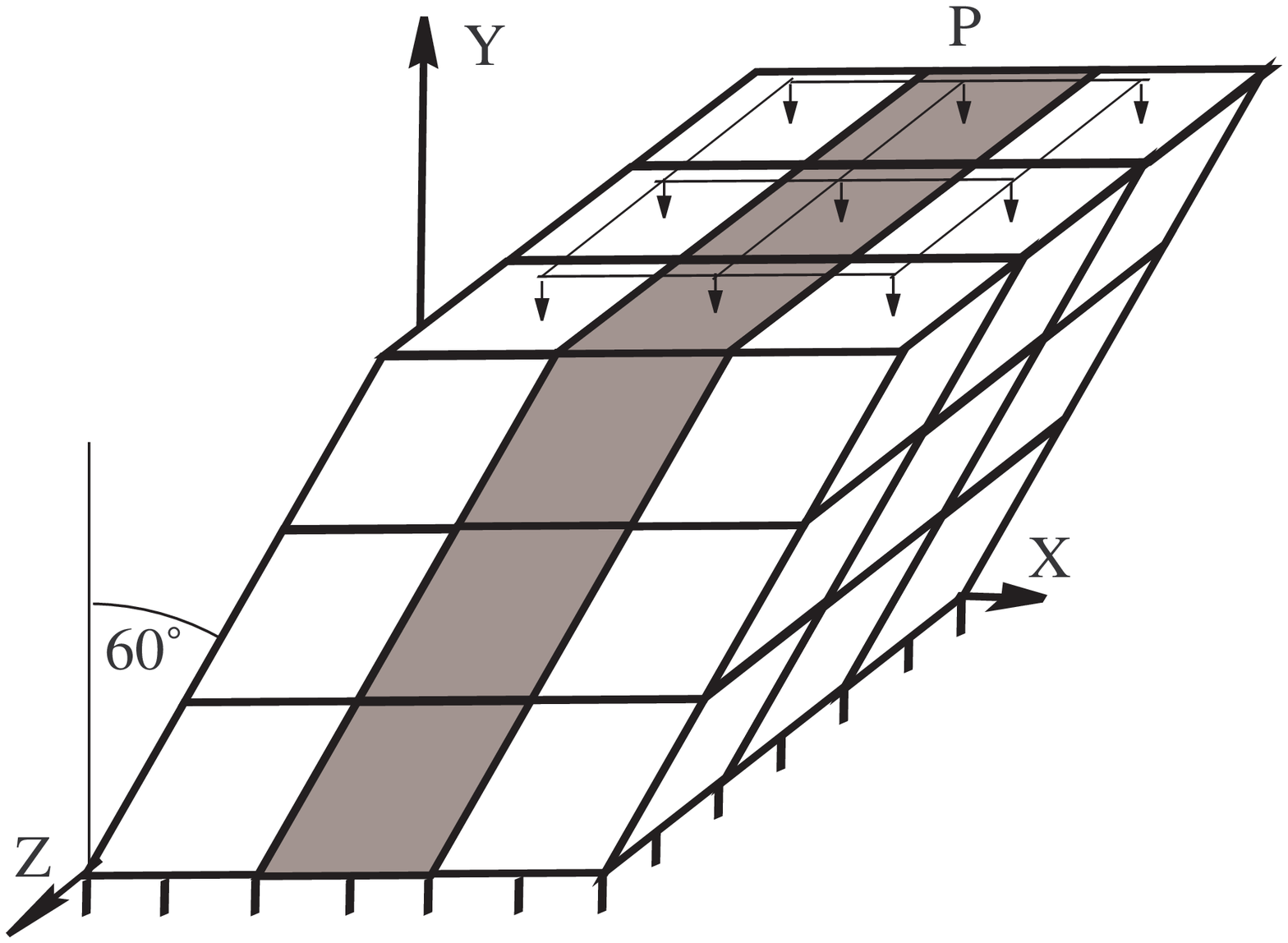}
\caption{Slanted cube with heterogeneous layers}\label{fig:cubeh360}
\end{minipage}
\end{figure}

\paragraph{Cube and slanted cube with heterogeneous layers: }
We now assess the new initialization strategy
on two other configurations of the cube in order to evaluate
the influence of the geometry and of the repartition of heterogeneities.
The structures depicted in
Figure \ref{fig:cubeh3} and \ref{fig:cubeh360} are similar to
 fig. \ref{fig:cube} but with different material distribution. Also, for the problem
described in
Figure \ref{fig:cubeh360}, the cube has been slanted by 60 degrees.
In
table \ref{tab:prob2-3} we report the number of iterations
when using the BDD solver and when applying FETI with the standard and the new initialization.

\begin{table}[ht]\centering
\begin{tabular}{|c|c|c|c|}\hline \multicolumn{2}{|c|}{Approach} & Cube, fig. \ref{fig:cubeh3}
&  Slanted cube, fig. \ref{fig:cubeh360}  \\
\multicolumn{2}{|c|}{ } &  Num. of it. & Num. of it.   \\ \hline  \hline
\multicolumn{2}{|c|}{BDD} & $19$ & $73$ \\ \hline \hline
FETI       & No Splitting             & $20$  & $85$ \\ \cline{2-4}
Dirichlet  & Classical Splitting  & $21$ & $85$ \\ \cline{2-4}
$P(\tilde{F}^{-1}_{I_D})$         & New Initialization    & $19$ & $73$ \\ \hline  \hline
FETI        & No Splitting             & $22$ & $88$ \\ \cline{2-4}
Dirichlet  & Classical Splitting  & $22$ & $88$ \\ \cline{2-4}
$P(W) $       & New Initialization   & $22$ & $85$ \\ \hline  \hline
FETI        & No Splitting            & $92$ & $153$ \\ \cline{2-4}
Dirichlet  & Classical Splitting  & $92$ & $154$ \\ \cline{2-4}
$P(I)$          & New Initialization    & $90$ & $159$ \\ \hline
\end{tabular}\caption{Performance results on problem fig. \ref{fig:cubeh3} and \ref{fig:cubeh360}}\label{tab:prob2-3}
\end{table}

For the straight cube of fig.\  \ref{fig:cubeh3},
the use of the novel initialization leads to only small improvements.
Observing that the FETI solver with Dirichlet projector and standard splitting converges
nearly as fast as the BDD, it is clear that the new splitting strategy
can not have a significant effect.

For the slanted cube, due to the geometric distortion of the substructures and of the mesh,
the BDD convergence is better than when the FETI with the standard splitting is applied.
When the new initialization strategy is used, the convergence history of FETI is again
very similar to the convergence of the primal BDD.

\begin{figure}[ht]\centering
\psfrag{A}[B]{rubber}
\psfrag{B}{metal}
\includegraphics[width=0.8\textwidth,clip=true]{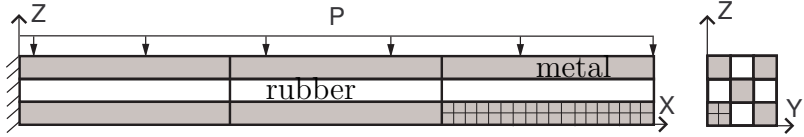}
\caption{$3 \times 3 \times 3$ substructures composite beam}\label{fig:beam}
\end{figure}

\paragraph{Nonlinear flexion of a composite beam}
As a last example, we analyze a non-linear problem solved through a sequence
of linearized systems.
The structure (fig.\ \ref{fig:beam}) is a slender beam of aspect ratio 9 with square
cross-section. It is made of longitudinal strips of
metal and rubber. The loading corresponds to an imposed pressure on one side of the beam.
Due to the presence of elastomer parts, the structure undergoes large deformations.
We chose a Kirchoff Saint-Venant model for the metal,
the characteristic coefficients of which are its Young modulus $E=20 000$ MPa
and Poisson's coefficient $\nu=0.3$.
A NeoHookian model is assumed for the rubber
characterized by a shear modulus $G=2.0$ MPa and a compressibility modulus $K=2000$ MPa.
The structure is decomposed into $27$ monomaterial paralellepipedic subdomains and
each substructure is meshed in $2\times2\times18$ cubic elements.
In order to handle the quasi-incompressibility of the elastomer, a mixed finite
element formulation is considered where the
pressure field and displacements are discretized independently. We choose
a $Q_2-P_1$ hexaedral element with 27 displacement nodes and 4 internal pressure nodes.
Details on this formulation and on the practicalities
for applying it properly in simulation  can be found in \cite{GOSSELET:2002:DDM}.

The complete model has $55 300$ degrees of freedom of which $16 400$
belong to the interface.
Due to the behaviour of the elastomer, the problem is highly nonlinear.
Newton-Raphson iterations are performed where
the tangent matrix is updated at every step.
The pressure loading is $5$ bars. The stopping criterion is set to
$10^{-3}$ for the relative primal residual of FETI when solving
the linearized systems. For the Newton-Raphson iterations, the tolerance
for the relative residual of the nonlinear equations is set to $10^{-4}$,
so that $5$ Newton-Raphson iterations and thus 5 linear solves
must be performed.

\begin{figure}[ht]\centering
\psfrag{I}[][][2.]{number of iterations}
\psfrag{N}[][][2.]{System number}
\psfrag{E}[Tr][][1.5]{FETI \small{$P(\tilde{\Fbold}^{-1}_{I_D})$}}
\psfrag{F}[Tr][][1.5]{improved initialization FETI \small{$P(\tilde{\Fbold}^{-1}_{I_D})$}}
\psfrag{V}[Tr][][1.5]{FETI \small{$\P(\Wbold)$}}
\psfrag{W}[Tr][][1.5]{improved initialization FETI \small{$\P(\Wbold)$}}
\psfrag{P}[Tr][][1.5]{BDD}
\includegraphics[angle=-90,width=0.8\textwidth]{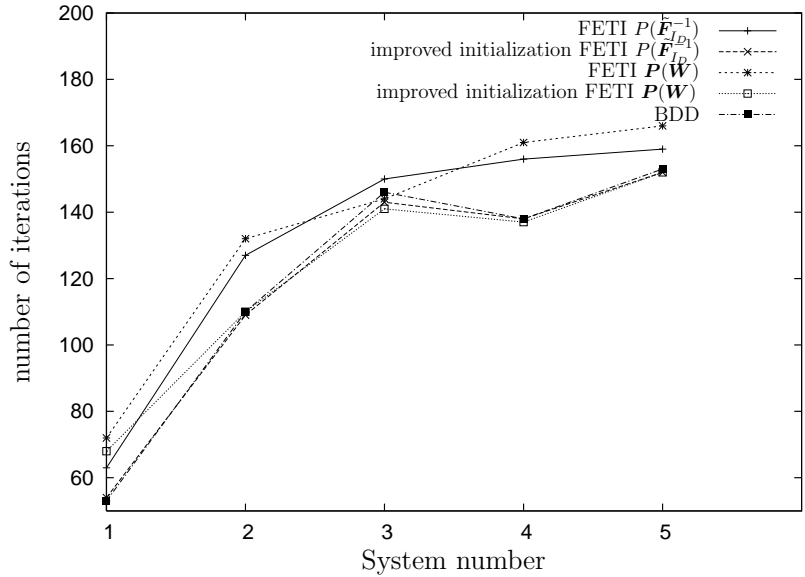}
\caption{Performance results for the flexion of the composite beam}\label{fig:iternl}
\end{figure}

\noindent Figure \ref{fig:iternl} indicates the number of
iterations required when solving the linearized systems by the
classical FETI method, by FETI with improved initialization and by
the primal BDD. The Dirichlet and diagonal projectors are applied,
and the Dirichlet preconditioner with stiffness scaling is used.
The significant increase of the number of iterations between the
first and the following linearized systems can be explained by the
loss of positivity of the tangent matrix which is mostly due to
incompressibility. Indeed, the conjugate gradient algorithm (with
full re-orthoganilization) applied on the interface problem in the
FETI and in the BDD methods remains applicable for non positive
matrices but its convergence is significantly slowed down
\cite{PAIGE:1995:ASE}. \footnote{Often, non-positivity is due to
the mutation of former null-modes (rotations) to negative modes,
it can be handled by the introduction of these modes as
constraints in Krylov-augmented algorithms
\cite{FARHAT:2000:FETINL}. In our case, non-positivity is mostly
due to the behaviour of the rubber and the strategy described
above is non-relevant. An efficient strategy based on the
approximation of negative modes can be found in
\cite{GOSSELET:2002:SRK}.} As observed from figure
\ref{fig:iternl}, the new initialization enables FETI to achieve
performance results which are very similar to the BDD method.
Comparing the total number of iterations, the classical FETI
approach requires $10\%$ more iterations than BDD while the FETI
method with new initialization requires slightly less iterations
than the BDD. From figure \ref{fig:iternl} we also observe that
the new initialization technique improves the performance of FETI
both for the Dirichlet and the diagonal projector. For this
particular structure with regular geometry, the cost effective
diagonal projector yields a convergence rate very similar to the
convergence rate obtained with the more computationally intensive
Dirichlet projector, except for the very first linearized system
solve.

%Comparing the total CPU time, the classical FETI method is
%$5\%$ slower than the BDD, while FETI with new initialization is $3\%$ faster.
Our numerical experiments indicate that for a large class of nonlinear problems
such as the one depicted here,
the new initialization leads to a small but non-negligible gain in terms of
number of iterations and CPU time.
Another important beneficial effect of the new starting procedure for FETI
comes from the fact that, since the initial residual of the iterations on the interface problem
are several orders of magnitude lower with the new initialization, stagnation of the residual
of the FETI iterations which often happens when dealing with higher nonlinearity (higher loading)
is significantly delayed, so that in practice restarting of the iteration can be avoided.
%We also observed that when willing to achieve better accuracy (by decreasing the convergence criteria)
%or when dealing with higher nonlinearity (higher loading), the new initialization makes FETI much more robust:
%indeed for such problems, stagnation often occurs with classical FETI requiring various restarting procedures,
% while new initialization avoids most of the numerical difficulties achieving then much better performance results.

%

%%%%%%%%%%%%%%%%%%%%%%%%%%%%%%%%%%%%%%%%%%%%%%%%%%%%%%%%%%%%%%%%%%%%%
%\section{Similarity and differences between primal and dual Schur complement procedures}
%
%options
%k scaling
%descent on different problems

%%%%%%%%%%%%%%%%%%%%%%%%%%%%%%%%%%%%%%%%%%%%%%%%%%%%%%%%%%%%%%%%%%%%%
\section{Conclusions}

For some particular structural problems such as those exhibiting strong heterogeneities
and geometric distortion, the Finite Element Tearing and Interconnecting (FETI) solver
can yield poor convergence compared to the conceptually similar Balanced Domain Decomposition
(BDD)
solver. In those cases, the bad performance of FETI can be traced back to high
initial residual in the iterations on the interface problem. The initial residual is
strongly related to the way the applied forces are split on the subdomain interface boundaries.

In this paper, we propose a novel strategy to split the applied forces between the subdomains.
We propose to split the statically condensed interface force according to
interface diagonal stiffness. This leads to
building a more efficient initial estimate for the
interface connecting forces.
The new initialization for the FETI iterations mainly involves
computing statically condensed forces on the interface and can thus be performed
at a computational cost equivalent to less than half the cost of a full FETI iteration.

The numerical examples described in this paper indicate that for the problems where
the primal BDD method outperforms FETI due to unexpected high
initial residual, the new initialization strategy builds a better starting estimate of
the interface forces and, in turn, to an initial residual similar to the BDD residual.
The FETI method then converges in a manner very similar to the BDD method.

The novel starting strategy never deteriorates the FETI convergence and leads to significant
improvements in some pathological cases. Therefore we suggest
to use the presented initialization as default in FETI solvers.

With the proposed initialization for FETI, the FETI method and the BDD solver lead to
similar convergence and computational costs for complex problem where the Dirichlet
projector is required. For problems where the simplified FETI preconditioners and
projectors can be used without significantly deteriorating the convergence of the
interface iterations, FETI is often found to be more efficient in terms of
overall computing cost.

%%%%%%%%%%%%%%%%%%%%%%%%%%%%%%%%%%%%%%%%%%%%%%%%%%%%%%%%%%%%%%%%%%%%%

\section*{Acknowledgements}
Part of this work was performed while the third author was visiting the LM2S
with the financial support of the Ecole Normal Sup\'erieure de Cachan, France. The first two authors
acknowlegde computational resources support from the Centre Informatique National Enseignement Sup\'{e}rieur (CINES) and the
P\^ole de Calcul Paris Sud.

\bibliography{\the\files}

\begin{thebibliography}{10}
\expandafter\ifx\csname url\endcsname\relax
  \def\url#1{\texttt{#1}}\fi
\expandafter\ifx\csname urlprefix\endcsname\relax\def\urlprefix{URL }\fi

\bibitem{MANDEL:1993:BAL}
J.~Mandel, Balancing domain decomposition, Comm. Appl. Num. Meth. Engrg. 9
  (1993) 233--241.

\bibitem{LETALLEC:1994:DDM}
P.~L. {T}allec, Domain-decomposition methods in computational mechanics,
  Computational Mechanics Advances 1~(2) (1994) 121--220, north-Holland.

\bibitem{FARHAT:1991:FETI}
C.~Farhat, F.-X. Roux, A method of finite tearing and interconnecting and its
  parallel solution algorithm, International J. Numer. Methods Engineering 32
  (1991) 1205--1227.

\bibitem{FARHAT:1994:ADV}
C.~Farhat, F.~X. Roux, Implicit parallel processing in structural mechanics,
  Computational Mechanics Advances 2~(1) (1994) 1--124, north-Holland.

\bibitem{RIXEN:2001:PARA}
D.~Rixen, Encyclopedia of Vibration, Academic Press, 2002, Ch. Parallel
  Computation, pp. 990--1001, iSBN 0-12-227085-1.

\bibitem{MANDEL:1996:OPT}
J.~Mandel, R.~Tezaur, Convergence of a substructuring method with {L}agrange
  multipliers, Numerische Mathematik 73 (1996) 473--487.

\bibitem{KLAWONN:2001:BAB}
A.~Klawonn, O.~Widlund, {FETI} and {N}eumann-{N}eumann iterative substructuring
  methods: Connections and new results, Comm. Pure App. Math. 54~(1) (2001)
  57--90.

\bibitem{FRAGAKIS:2002:FETIP1}
Y.~Fragakis, M.~Papadrakakis, A unified framework for formulating domain
  decomposition methods in structural mechanics, Tech. rep., Institute for
  Structural Analysis \& Seismic Research, Athens, Greece (2002).

\bibitem{PARK:1997:ALG}
K.~Park, M.~Justino, C.~Felippa, An algebraically partitioned {FETI} method for
  parallel structural analysis: Algorithm description, International J. Numer.
  Methods Engineering 40~(15) (1997) 2717--2737.

\bibitem{RIXEN:1999:AFETI}
D.~Rixen, C.~Farhat, R.~Tezaur, J.~Mandel, Theoretical comparison of the feti
  and algebraically partitioned feti methods, and performance comparisons with
  a direct sparse solver, International J. Numer. Methods Engineering 46~(4)
  (1999) 501--534.

\bibitem{RIXEN:2002:FETI1MPC}
D.~Rixen, Extended preconditioners for {FETI} method applied to constrained
  problems, Internat. J. Num. Meth. Engin. 54~(1) (2002) 1--26.

\bibitem{LETALLEC:1991:DDM}
P.~L. {T}allec, Y.-H.~D. {R}oeck, M.~Vidrascu, Domain-decomposition methods for
  large linearly elliptic three dimensional problems, J. of Computational and
  Applied Mathematics 34 (1991) 93--117, elsevier Science Publishers,
  Amsterdam.

\bibitem{FARHAT:2001:FETI_DP}
C.~Farhat, M.~Lesoinne, P.~LeTallec, K.~Pierson, D.~Rixen, {FETI-DP}: a
  dual-primal unified {FETI} method - part i: a faster alternative to the
  two-level {FETI} method, International J. Numer. Methods Engineering 50~(7)
  (2001) 1523--1544.

\bibitem{DELABOURDONNAYE:1998:FETIH}
A.~de~La~Bourdonnaye, C.~Farhat, A.~Macedo, F.~Magoules, F.-X. Roux, Advances
  in Computational Mechanics with High Performance Computing, Civil-Comp Press,
  Edinburgh, United Kingdom, 1998, Ch. A method of finite element tearing and
  interconnecting for the Helmholtz problem, pp. 41--54.

\bibitem{RIXEN:1998:SUPERL}
D.~Rixen, C.~Farhat, A simple and efficient extension of a class of
  substructure based preconditioners to heterogeneous structural mechanics
  problems, Internat. J. Num. Meth. Engin. 44~(4) (1999) 489--516.

\bibitem{Bhardwaj:1999:PQ_LRV}
M.~Bhardwaj, D.~Day, C.~Farhat, M.~Lesoinne, K.~Pierson, D.~Rixen, Application
  of the {FETI} method to {ASCI} problems: Scalability results on a
  thousand-processor and discussion of highly heterogeneous problems,
  International J. Numer. Methods Engineering 47~(1-3) (2000) 513--536.

\bibitem{GOSSELET:2002:DDM}
P.~Gosselet, C.~Rey, P.~Dasset, F.~Lene, A domain decomposition method for
  quasi incompressible formulations with discontinuous pressure field, Revue
  Européenne des Elements Finis 11 (2002) 363--377.

\bibitem{PAIGE:1995:ASE}
C.~Paige, B.~Parlett, H.~van~der Vorst, Approximate solutions and eigenvalue
  bounds from krylov subspaces, Numerical Linear Algebra with Applications
  2~(2) (1995) 115--133.

\bibitem{FARHAT:2000:FETINL}
C.~Farhat, K.~Pierson, M.~Lesoinne, The second generation {FETI} methods and
  their application to the parallel solution of large-scale linear and
  geometrically non-linear structural analysis problems, Comput. Meth. Appl.
  Mech. Engin. 184~(2-4) (2000) 333--374.

\bibitem{GOSSELET:2002:SRK}
P.~Gosselet, C.~Rey, On a selective reuse of krylov subspaces in newton-krylov
  approaches for nonlinear elasticity, in: Proceedings of the 14th Conference
  on Domain Decomposition Methods, Mexico 2002.

\end{thebibliography}

\end{document}